\documentclass[11pt]{amsart}

\usepackage[all]{xy}
\usepackage{xcolor}

\numberwithin{equation}{section}

\newtheorem{theorem}[equation]{Theorem}

\newtheorem{cor}[equation]{Corollary}

\theoremstyle{definition}

\newtheorem{example}[equation]{Example}

\theoremstyle{remark}
\newtheorem{remark}[equation]{Remark}

\newcommand{\C}{\mathbb{C}}
\newcommand{\g}{\gamma}
\newcommand{\hg}{\widehat{\g}}
\newcommand{\h}{\lambda}
\newcommand{\hh}{\widehat{\h}}
\newcommand{\N}{\mathbb{N}}
\newcommand{\Z}{\mathbb{Z}}

\newcommand{\om}{\omega}

\newcommand{\hG}{\widehat{\Gamma}}
\newcommand{\hD}{\widehat{\Lambda}}

\newcommand\Ext{\operatorname{Ext}}
\newcommand\Homol{\operatorname{H}}
\newcommand\coh{\Homol}
\newcommand\ot{\otimes}

\newcommand\ott{\otimes^{t}\!}
\newcommand{\ttp}{A\ot^{\tau}\! B}
\newcommand{\ttpt}{A\ott B} 

\newcommand\Hom{\operatorname{Hom}}

\DeclareMathOperator{\Ker}{Ker}

\DeclareMathOperator{\Tor}{Tor}

\newcommand{\DOT}{\setlength{\unitlength}{1pt}
     \begin{picture}(2.5,2)(1,1)\put(2,3){\circle*{2}}\end{picture}}
\newcommand{\bu}{\DOT}

\begin{document}

\thanks{X.~W.~was partially supported by the Simons Foundation Travel Support for
Mathematicians award SFI-MPS-TSM-00025285, and  
S.~W.~by award~SFI-MPS-TSM-00013813.}

\title[Cohomology of twisted tensor products]
{Cohomology of twisted tensor products\\ over commutative rings}

\author{Xingting Wang}
\address{Department of Mathematics\\Louisiana State University\\Baton Rouge, LA 70803}
\email{xingtingwang@lsu.edu}

\author{Sarah Witherspoon}
\address{Department of Mathematics\\Texas A\&M University\\College Station, TX 77843}
\email{sjw@tamu.edu}

\date{19 November 2026}

\begin{abstract}
Over a commutative ground ring,
cohomology of a twisted tensor product of two rings
mixes their cohomology via the twisting.
We focus on bicharacter twistings and give explicit descriptions
when the ground ring is itself a module. 
We state our main results for 
ground rings that are hereditary and contain
roots of unity corresponding to finite grading groups.
As a class of examples, 
we compute cohomology rings of Taft orders over cyclotomic integers. 
\end{abstract}

\maketitle

\section{Introduction}

Cohomology of a finite dimensional Hopf algebra (order), with coefficients 
in a commutative ground ring $R$, is a graded commutative ring.
In specific settings, much is known about this cohomology, 
from early finite generation results for finite group rings
over noetherian rings~\cite{Evens} to recent
results for finite group schemes~\cite{GP,vdK}.
Other recent related work with commutative ground rings $R$
includes~\cite{Be,BL,Magland}.
By contrast, for many types of finite dimensional Hopf algebras,
finite generation results have mostly focused on the case that 
$R$ is a field.
Sometimes they include detailed structural results on the cohomology
rings in addition, such as in~\cite{BNPP14,GK93}.
It is valuable to extend knowledge to other ground rings
that are natural choices in some settings. 
As well, working over an integral domain
for example may allow comparison of
cohomology over a field of fractions
with cohomology over a residue field.
Such comparisons can be useful 
in proving cohomology is finitely generated; 
see, for example,~\cite{KV} on small quantum groups
and general results in~\cite{NWW}. 

This note contributes to understanding 
cohomology of Hopf algebras (orders)
over commutative ground rings other than fields.
We adapt some results of~\cite{BW} from Hochschild cohomology
over a field to augmented algebra (order) cohomology over a 
more general commutative ground ring $R$.
Specifically, we describe the cohomology of a twisted tensor product ring,
where the twisting is given by a bicharacter on grading groups,
in terms of cohomology groups of the two factor rings.
Our main results require the ground ring $R$ to be hereditary 
(used in applying the K\"unneth Theorem) and to contain roots of unity 
of order corresponding to orders of finite grading groups.
We describe the
graded $R$-module structure of cohomology in Theorem~\ref{thm:main}
in case $R$ is a graded module for one of the two constituent rings.
This description mixes graded components and twisted coefficients.  
Under the additional hypothesis that $R$ is a graded module
for both constituent rings, we give in Theorem~\ref{thm:cup}
the cohomology of a twisted tensor product ring 
as a twisted tensor product of the cohomology rings of the constituent rings.
We are however also interested in examples
that do not satisfy this additional hypothesis of Theorem~\ref{thm:cup}.
Section~\ref{sec:Taft} is devoted to such a class of examples,
the Taft orders over cyclotomic integers.
There we use explicit chain level constructions to give the ring structure
of cohomology of Taft orders in Theorem~\ref{thm:Taft-cup}. 
In particular, these cohomology rings are {\em not} twisted tensor products
of cohomology rings of constituent rings. 
We view these Taft orders as prototypical of 
a much larger class of examples to which our techniques will apply, the
Hopf orders arising as bosonizations of Nichols orders.
These include in particular the analogs over $R$ of quantum Borel
subalgebras of small quantum groups $u_q(\mathfrak{g})$. 

Our results and examples should be compared with~\cite[Theorem 5.5]{NWW},
which is missing a hypothesis.
In~\cite{NWW-corrigendum} it is shown that under an additional assumption,
namely the vanishing of certain torsion $R$-modules,
the cohomological finite generation result of~\cite[Theorem 5.5]{NWW} is correct. 
By contrast, we show in~Section~\ref{sec:Taft}
that cohomology of a Taft order has some torsion elements. 
Nonetheless, as we also show there, there are other 
techniques available to conclude cohomological finite generation 
over cyclotomic integers. 

We begin by recalling twisted tensor products in Section~\ref{sec:ttp}.
In Theorem~\ref{thm:ttpr} we generalize,
from fields to some commutative ground rings,  
a construction from~\cite{Shepler-W}
of a twisted tensor product projective resolution.
Section~\ref{sec:additive} is on the graded $R$-module structure of cohomology
of twisted tensor products, Section~\ref{sec:cup} is on cup products
under additional hypotheses, and Section~\ref{sec:Taft} contains the Taft order examples.

Throughout, $R$ will be a commutative ring, and for the main results
we will impose additional constraints on $R$. 
For notational convenience, we write $\ot$ in place of $\ot_R$
for maps and elements; it will be clear from context that
the tensor product is over $R$. 


\section{Twisted tensor products and resolutions}\label{sec:ttp}

Let $R$ be a commutative ring.
We will work with an {\em $R$-order} $A$ that is a ring containing $R$ 
as a subring of its center, under which containment $A$ is a finitely generated
free $R$-module.
We say that $A$ is {\em augmented} if there is an $R$-order homomorphism
$\varepsilon_A: A \rightarrow R$ that takes the subring $R$ of $A$ identically to itself. 
Let $A$ and $B$ be such $R$-orders. 
Denote their multiplication maps by
$m_A: A\ot_R A\rightarrow A$ and $m_B: B\ot_R B\rightarrow B$.

A {\em twisting map} $\tau: B\ot_R A\rightarrow A\ot_R B$ 
is a bijective $R$-module homomorphism for which there is an equality of
$R$-module homomorphisms from $B\ot_R B\ot_R A\ot_R A$ to 
$A\ot_R B$,
\begin{equation}\label{eqn:twisting-mult}
\tau \circ (m_B \ot m_A) = 
(m_A \ot m_B) \circ (1 \ot \tau \ot 1)
\circ (\tau \otimes \tau) \circ (1 \ot \tau \ot 1) , 
\end{equation}
and $\tau(1_B\ot a) = a\ot 1_B$ and 
$\tau(b \ot 1_A) = 1_A \ot b$ for all $a\in A$ and $b\in B$. 
(Here, $1_A$ and $1_B$ denote the multiplicative identities 
of $A$ and $B$, respectively, and $1$ generically denotes an identity map.) 
Define the 
{\it twisted tensor product $R$-order}
$\ttp$ to be the free $R$-module $A\otimes_R B$ 
with multiplication $m_\tau: (\ttp)\ot_R (\ttp) \rightarrow (\ttp)$ 
defined by $m_\tau =  (m_A\ot m_B) (1\ot \tau\ot 1)$. 
The defining properties of the twisting map $\tau$ ensure that $\ttp$ 
is itself an $R$-order with $R\cong R\ot_R R$ a central subring. 

We will construct a twisted tensor product of resolutions,
adapting results of~\cite{Shepler-W} from fields to more general
commutative ground rings $R$.
We will assume further for the main theorems
that $R$ is {\em hereditary}, that is, $R$-submodules of
projective $R$-modules are projective. 
All $A$-modules and $B$-modules in the constructions below are 
assumed to be free as $R$-modules. 
Later, we will encounter torsion modules when taking cohomology.

Let $M$ be an $A$-module with module action given by  
$\rho_{A,M}: A\ot_R M\rightarrow M$.
We say that $M$ is 
{\em compatible with the twisting map $\tau$} if there is a bijective
$R$-module homomorphism $\tau^{}_{B,M}: B\ot_R M\rightarrow M\ot_R B$ 
for which there are equalities of 
$R$-module homomorphisms on $B\ot_R B\ot_R M$ and on $B\ot_R A\ot_R M$: 
\begin{eqnarray*}
  \tau^{}_{B,M}(m_B\ot 1) & = & (1\ot  m_B)(\tau^{}_{B,M}\ot 1)(1\ot \tau^{}_{B,M}) ,\\
  \tau^{}_{B,M}(1\ot \rho^{}_{A,M})&  =&
  (\rho^{}_{A,M}\ot 1)(1\ot \tau^{}_{B,M})  (\tau \ot 1) . 
\end{eqnarray*}
If $M$ is compatible with $\tau$ and $N$ is any $B$-module,
then the tensor product $M\ot_R N$ is an $\ttp$-module by applying 
$1\ot\tau_{B,M}\ot 1$ to $\ttp\ot_R M\ot_R N$ followed by the 
$A$- and $B$-module actions.

We are particularly interested in $R$ itself as a module: 
If $A$ and $B$ are augmented, 
give $R$ the structure of an $A$-module and of a $B$-module via 
their augmentation maps $\varepsilon_A$ and $\varepsilon_B$.
Assume the $A$-module $R$ is compatible with $\tau$.
It follows that the $R$-module $R\cong R\ot_R R$
is an $\ttp$-module under $\varepsilon_A\ot \varepsilon_B$.
In this way, $\ttp$ is also an augmented $R$-order.

\subsection*{Twisted tensor product resolutions}
Let $M$ be an $A$-module that is compatible with
the twisting map $\tau$. 
Let $P$ be a projective resolution of $M$ as an $A$-module:
\[
  P: \hspace{1cm} \cdots \longrightarrow P_2
     \longrightarrow P_1 \longrightarrow P_0 \longrightarrow M
     \longrightarrow 0 .
\]
We say that $P$ 
is {\em compatible with the twisting map $\tau$} if
each $A$-module $P_i$ in $P$ is compatible with $\tau$ 
via an $R$-module homomorphism $\tau^{}_{B,P_i}$
in such a way that the resulting map 
$\tau^{}_{B,P}: B\ot_R P\rightarrow P\ot_R B$ 
is an $R$-module chain map lifting
$\tau^{}_{B,M}: B\ot_R M \rightarrow M\ot_R B$.

We give the following generalization of results  
in~\cite[Section 5]{Shepler-W}.

\begin{theorem}\label{thm:ttpr}
Let $R$ be a hereditary commutative ring.
Let $P$ (respectively $Q$) be a projective resolution of 
an $A$-module $M$ (respectively $B$-module $N$).
If $M$ and $P$ are compatible with $\tau$
and each $P_i\ot_R Q_j$ is a projective $\ttp$-module, 
then the total complex of $P\ot_R\! Q$ is a projective
$\ttp$-resolution of $M\ot_R N$.
\end{theorem}

\begin{proof}
The complexes $P$ and $Q$ consist of free $R$-modules by assumption.
Their tensor product complex $P\ot_R Q$ again consists of free $R$-modules.
Apply the K\"unneth Theorem~\cite[Theorem 3.6.3]{Weibel};  
since $R$ is hereditary, each of the 
$R$-modules $d_i(P_i)$, as well as each $P_i$, 
is indeed flat (in fact projective) over $R$
(where $d_i$ denotes the differential from $P_i$ to $P_{i-1}$).
Thus for each $n$ there is an exact sequence
\[
   0\rightarrow \bigoplus_{i+j=n}
  (\coh_i(P)\ot_R \coh_j(Q)) \rightarrow
   \coh_n(P\ot_RQ) \rightarrow
   \bigoplus_{i+j=n-1}
   \Tor^R_1(\coh_i(P), \coh_j(Q)) \rightarrow 0 .
\]
The homology of the complex $P$ (respectively, $Q$) is
concentrated in degree~0, where it is $M$ (respectively, $N$).
Since $M$ and $N$ are free (hence flat) over $R$ by our standing
assumptions, $\Tor^R_1(M,N)=0$. Thus the above Tor term vanishes and
the homology of the total complex of $P\ot_{R} Q$ is also
concentrated in degree 0, where it is $M\ot_R N$. 
It follows that the augmented complex 
Tot$(P\ot_R Q) \rightarrow M\ot_R N\rightarrow 0$
is exact.
By hypothesis, it consists of projective $\ttp$-modules, and so
it is a projective resolution of $M\ot_R N$. 
\end{proof}

\begin{remark}
In practice, for many standard types of resolutions and twisting maps, 
the projectivity of each module $P_i\ot_R Q_j$ is straightforward to check. 
See~\cite[Theorem~5.12]{Shepler-W} for some sufficient conditions
under which such a twisted tensor product, over a field, is a
projective resolution.
These conditions generalize to commutative ground rings in a straightforward way,
but we will not need that here. 
Note that bar resolutions always satisfy the compatibility hypothesis
of Theorem~\ref{thm:ttpr}, and we will use these in 
Sections~\ref{sec:additive} and~\ref{sec:cup} below.  
See Section~\ref{sec:Taft} for examples that are not bar resolutions.
\end{remark}

\subsection*{Twisting by a bicharacter} 
For our main results in the next sections, we will work in the special setting of 
twistings given by bicharacters on grading groups, as
we recall next. 

Let $\Gamma$ and $\Lambda$ be finite abelian groups, written multiplicatively.
Let $R$ be a commutative ring containing primitive
roots of unity of orders $|\Gamma|$ and $|\Lambda|$.
Let $A$ and $B$ be $R$-orders that are graded by the groups
$\Gamma$ and $\Lambda$, that is, 
\[
  A = \bigoplus_{\g\in \Gamma}A^{\g} \ \ \ \ 
  \mbox{ and } \ \ \ \ B = \bigoplus_{\h\in\Lambda} B^{\h}
\]
as $R$-modules, with
$A^{\g}A^{\g '}\subseteq A^{\g \g '}$ and $B^{\h}B^{\h '} \subseteq B^{\h \h '}$
for all $\g , \g'\in \Gamma$ and $\h , \h'\in \Lambda$.
Assume the multiplicative identity $1_A$ of $A$ is contained
in the identity component $A^{1_{\Gamma}}$ of $A$, and similarly for $B$.

Denote by $R^{\times}$ the multiplicative group of units of $R$. 
Let $t: \Gamma\times\Lambda \rightarrow R^{\times}$ be a bicharacter, that is
\[
t(\g\g',\h) = t(\g,\h)t(\g',\h) \ \ \mbox{ and } \ \
t(\g,\h\h') = t(\g,\h)t(\g,\h')
\]
for all $\g,\g'\in\Gamma$ and $\h,\h'\in\Lambda$.
We will also write $t(a,b) = t(|a|,|b|)$ when $a\in A$ and $b\in B$
are homogeneous of degrees $|a|$ and $|b|$ in $\Gamma$ and $\Lambda$,
respectively.
The map $\tau: B\ot_R A\rightarrow A\ot_R B$ given on 
homogeneous elements $a\in A$ and $b\in B$ by 
\[
   \tau(b\ot a) = t(a,b) a\ot b
\]
is a twisting map.
Denote by $\ttpt$ the corresponding twisted tensor product $R$-order.
Specifically, multiplication on $\ttpt$ is given by 
\[
(a\ot b)(a'\ot b') = t(a',b) aa'\ot bb'
\]
for all elements $a,a'\in A$ and $b,b'\in B$,
with $a', b$ homogeneous. 


\section{Cohomology of bicharacter twisted tensor products}\label{sec:additive}

Our main results in this section are versions of some of those in~\cite{BW}
in a little bit different context.
In particular, here we replace 
Hochschild cohomology with augmented order cohomology
and a ground field with a more general ground ring.
We will use the twisted tensor product resolution of Theorem~\ref{thm:ttpr}
for the choices $M=R$, $N=R$. 

First we recall some relevant automorphisms from~\cite{BW}.
We assume $R$ contains primitive roots of unity
of orders $|\Gamma|$ and $|\Lambda|$ for finite abelian groups 
$\Gamma$ and $\Lambda$, and $A$ and $B$ are augmented $R$-orders
graded by $\Gamma$ and $\Lambda$, respectively.
Let
\[
\hG = \Hom_{{\rm{gp}}} (\Gamma , R^{\times}) ,
\]
the dual group of characters of $\Gamma$ with values in
the group $R^{\times}$ of units of $R$,
that is all group homomorphisms from $\Gamma$ to $R^{\times}$.
Similarly define $\hD$.
A bicharacter $t:\Gamma\times\Lambda\rightarrow R^{\times}$ induces group homomorphisms 
$\Lambda  \longrightarrow \hG$
and $\Gamma  \longrightarrow \hD$, given by
$\h \mapsto \hh$ and $\g   \mapsto \hg$ 
where 
\[
   \hh = t( - , \h) \ \ \ \mbox{ and } \ \ \ \hg = t(\g, - )
\] 
for all $\gamma\in\Gamma$ and $\lambda\in\Lambda$. 
These maps in turn induce group homomorphisms from $\Lambda$
and $\Gamma$ to the groups of $R$-order automorphisms of $A$ and $B$
as follows. 
We use the same notation for the corresponding $R$-order automorphisms:
\[
\begin{aligned}
  &\hh (a) = \hh (\g) a = t(\g,\h)a , \ \ \
          &\hg (b)  = \hg (\h)b = t(\g,\h)b \ \ \
\end{aligned}
\]
for all $a\in A^{\g}$, $\g\in \Gamma$,
$b\in B^{\h}$, and $\h\in\Lambda$.

Let $M$ be an $A$-module.
For each $\h\in\Lambda$,
denote by $M_{\hh}$ the $A$-module that is $M$ as an $R$-module
with action of $A$ twisted by the automorphism $\hh$:
\begin{equation}\label{eqn:mod-twist}
a\cdot_{\hh} m = \hh (a) m 
\end{equation}
for all $a\in A$, $m\in M$.
Similarly, if $N$ is a $B$-module and $\g\in\Gamma$,
denote by $N_{\hg}$ the $B$-module with action twisted by
the automorphism $\hg$ of $B$.

Assume that $\ker \varepsilon_A$ contains all homogeneous
elements in $A$ that are not in the identity component $A^{1_{\Gamma}}$.
Then $R$ is a $\Gamma$-graded $A$-module concentrated in degree $1_{\Gamma}$, 
under the action induced by $\varepsilon_A$.
We do {\em not} assume that $R$ is a $\Lambda$-graded $B$-module; 
this asymmetry in our assumptions on the augmentation maps 
$\varepsilon_A, \varepsilon_B$ in theorem statements below corresponds to
the choice of left (versus right) modules and the position of 
$A$ as the left (versus right) factor in $\ttpt$.
The Taft order examples in Section~\ref{sec:Taft} below 
satisfy these assumptions, as we explain there.
They are prototypical of a much larger class of examples which 
include finite rank Hopf orders
that are semidirect products (bosonizations) of Nichols orders with finite groups,
such as analogs over $R$ of quantum Borel subalgebras of small
quantum groups $u_q(\mathfrak{g})$.

Give $R$ the structure of an $A\ott B$-module 
via the augmentation map $\varepsilon_{A\ott B}=\varepsilon_A\ot\varepsilon_B$ 
as in Section~\ref{sec:ttp}. 
  Let $P$ (respectively $Q$) be the 
  bar resolution of $R$ as an
  $A$-module (respectively as a $B$-module).
  That is,
  \begin{equation}\label{eqn:P-bar}
  P: \hspace{1cm}
  \cdots \stackrel{d_3}{\longrightarrow} A\ot_R A\ot_R A
  \stackrel{d_2}{\longrightarrow} A\ot_R A
  \stackrel{d_1}{\longrightarrow} A
  \stackrel{\varepsilon_A}{\longrightarrow} R \longrightarrow 0 ,
  \end{equation}
with action of $A$ given by multiplication on the leftmost tensor
  factor in each degree. The differentials are defined by 
  \[
  d_i (a_0\ot\cdots\ot a_i) =  (-1)^i\varepsilon_A(a_i) a_0\ot\cdots\ot a_{i-1} +
\sum_{k=0}^{i-1} (-1)^k
  a_0\ot\cdots\ot a_ka_{k+1}\ot \cdots \ot a_i 
     \]
for all $i\geq 1$ and $a_0,\ldots, a_i \in A$.
This resolution $P$ is compatible with $\tau$,
as may be seen by iterating the twisting map on $A$, $B$.
Each $P_i\ot_R Q_j$ is free as an $\ttpt$-module since
iterating the twisting gives 
isomorphisms of $R$-modules, 
  \begin{equation}\label{eqn:pqab}
  P_i\ot_R Q_j = A^{\ot_R (i+1)}\ot_R B^{\ot_R (j+1)} \cong
  (A\ott B) \ot_R (A^{\ot_R i}\ot_R B^{\ot_R j}) .
  \end{equation}
Thus by Theorem~\ref{thm:ttpr}, provided $R$ is hereditary, 
the total complex of $P\ot_R Q$ is a free $A\ott B$-resolution of $R$,
denoted here by $P\ott Q$.

The terms in the bar resolution $P$ of $R$ as
an $A$-module are $\Gamma$-graded with $\Gamma$-degrees 
\[
    |a_1\ot\cdots\ot a_n| = |a_1|\cdots |a_n| 
\]
for all homogeneous $a_1,\ldots, a_n\in A$. 
The $\Gamma$-grading on $A$ then yields a $\Gamma$-grading
on the cohomology $\coh(A,R) = \Ext^*_A (R,R) = \oplus_{i\geq 0}\Ext^i_A(R,R)$,
\[
\coh(A,R) \cong \bigoplus_{\g\in\Gamma}  \coh(A,R)^{\g} ,
\]
where the grading is induced by that at the chain level:
For $\gamma\in\Gamma$, set 
\begin{equation*}
\begin{aligned}
 & \Hom_R (A^{\ot n}, R)^{\gamma} = 
    \{ f\in \Hom_R(A^{\ot n}, R)\mid
    |f(a_1\ot\cdots\ot a_n)| = \gamma^{-1}  |a_1|\cdots |a_n| \\
  & \hspace{7.5cm}   \mbox{ for all homogeneous } a_1,\ldots,a_n\in A\}\\
  & \hspace{2.8cm} = \{ f\in \Hom_R(A^{\ot n},R)\mid
       f(a_1\ot\cdots\ot a_n)\neq 0 \mbox{ only if }
         |a_1|\cdots |a_n| = \gamma \} ,
\end{aligned}
\end{equation*}
the second equality holding since $R$ has degree~$1_{\Gamma}$ as an $A$-module.

To prove the following theorem, we use 
the twisted tensor product resolution with components~(\ref{eqn:pqab}) and 
adapt the proof of~\cite[Theorem 3.1]{BW}, 
making the changes from Hochschild cohomology for algebras over a field
to augmented order cohomology for orders over a hereditary commutative ring. 
The $B$-modules $R_{\hg}$ (see~(\ref{eqn:mod-twist}))
in the theorem statement below arise naturally,
for example there is an $R$-module homomorphism 
\[
  \Hom_A(A,R)^{\gamma} \ot_R \Hom_B(B, R_{\hg}) \stackrel{\sim}{\longrightarrow}
   \Hom_{\ttpt} (\ttpt , R)
\]
given by setting
$(f\ot g)(a\ot b) = f(a)g(b)$ for all $a\in A$, $b\in B$,
$f\in\Hom_A(A,R)^{\gamma}$, $g\in \Hom_B(B,R_{\hg})$,
as we explain next.
To see that this function $f\ot g$ is indeed an $\ttpt$-module
homomorphism, note that 
\begin{eqnarray*}
  (f\ot g) ((a'\ot b')(a\ot b)) & = & t(a,b') f(a'a) g(b'b) , \\
  (a'\ot b') ( (f\ot g) (a\ot b)) & = & (a'\ot b') (f(a)\ot g(b)) \\
 & = & a'f(a) \ot \hg(b') g(b) \ \ = \ \ f(a'a) \ot t(a,b') g(b'b)
\end{eqnarray*}
for homogeneous $a,a'\in A$ and $b,b'\in B$ since $f(a)\neq 0$ only if $|a|=\gamma$.
In comparison with~\cite{BO}, note that the automorphism $\hg$ is intrinsic in the isomorphism
of~\cite[Lemma 3.2]{BO},  
\[
  M\langle \g \rangle \ot^t \! N\langle \lambda \rangle
  \rightarrow (M\ot^t\! N)\langle \g , \lambda \rangle ,
\]
given by $m\ot n \mapsto  t(\g,n) m\ot n$.
Their factor $t(\g,n)$ is a value of $\hg$.

\begin{theorem}\label{thm:main}
  Let $R$ be a hereditary commutative ring that is a 
$\Gamma$-graded $A$-module concentrated in degree~$1_{\Gamma}$,
  and assume that for each $i\geq 0$, 
  $\coh^i(A,R)$ is flat over $R$.
  Then for each $m\geq 0$, 
  there is an isomorphism of $R$-modules:
  \[
  \coh^m(A\ott B , R)\cong \bigoplus_{i+j=m, \ \g\in\Gamma}
   \left( \coh^i(A, R)^{\g} \ot_R \coh^j(B, R_{\hg}) \right) .
  \]
\end{theorem}

\begin{proof}
(cf.~\cite[Theorem 3.1]{BW})
Let $P$ (respectively, $Q$) be the bar resolution of $R$
as an $A$-module (respectively, $B$-module); see e.g.~(\ref{eqn:P-bar}).
By Theorem~\ref{thm:ttpr} and the isomorphisms~(\ref{eqn:pqab}), $P\ott Q$ is a
projective $A\ott B$-module resolution of $R$.
Thus $\coh(A\ott B, R)$ is the
  cohomology of the cochain complex 
\[
   C_t(A,B,R) = \Hom_{A\ott B}(P\ott Q , R) .
\]
  Invoking the isomorphisms~(\ref{eqn:pqab}),  
  in each degree $m$, the $R$-module $C^m_t(A,B,R)$ is a direct sum, over
  $i+j=m$, of components
  \begin{equation}\label{eqn:homab}
  \Hom_{A\ott B}(A^{\ot_R  (i+1)}\ot_R  B^{\ot_R (j+1)}, R)
  \cong  \Hom_{R}(A^{\ot_R  i}\ot_R  B^{\ot_R j}, R) .
  \end{equation}
  We will further analyze this cochain complex.
    Let $C(A,R) = \Hom_A (P,R)$, which by hypothesis is graded by $\Gamma$.
    There are isomorphisms
    \[
    C^i(A,R) \cong  \bigoplus_{\g\in\Gamma}\Hom_R(A^{\ot_R i} , R)^{\g} . 
    \]
Next, for each $\g\in \Gamma$, let $C(B,R_{\hg}) = \Hom_B(Q, R_{\hg})$.
There are isomorphisms
  \[
   C^j(B,R_{\hg})  \cong  \Hom_R (B^{\ot_R j},R_{\hg}) .
   \]
   We will show that there is an isomorphism of cochain complexes of $R$-modules, 
   \begin{equation}\label{eqn:cacb}
   C_t(A,B,R) \cong  
  \bigoplus_{\g\in\Gamma} \left(C(A,R)^{\g} \ot_R C(B, R_{\hg})\right) .
  \end{equation}

For each $\g\in\Gamma$, define an $R$-module homomorphism
  \[
    \phi_{\g}:    C(A,R)^{\g} \ot_R C(B, R_{\hg}) \rightarrow C_t(A,B,R)
    \]
as follows. 
Let $f\in C^i(A,R)^{\g}$, $g\in C^j(B,R_{\hg})$, and set
    \[
    \phi_{\g}(f\ot g) (a_1\ot\cdots\ot a_i \ot  b_1\ot \cdots\ot b_j)
    = (-1)^{ij}  f(a_1\ot\cdots \ot a_i)g(b_1\ot\cdots\ot b_j) .
    \]
Now let $\phi = \sum_{\g\in\Gamma}\phi_{\g}$, to obtain a map
\[
  \phi:  \bigoplus_{\g\in\Gamma} \left(C(A,R)^{\g} \ot_R C(B, R_{\hg})\right)
  \rightarrow C_t(A,B,R) .
  \]
The domain and codomain of $\phi$ are $\N\times \Gamma$-graded $R$-modules. 
Since $A$ and $B$ are assumed to be 
free $R$-modules of finite rank and $R$ is hereditary, 
all graded components above are projective $R$-modules of finite rank.
It is now clear from its definition
that $\phi$ is an isomorphism of $\N\times\Gamma$-graded $R$-modules.
    We show next that $\phi$ is a chain map.
    As a consequence, $\phi$ will induce an isomorphism on cohomology.
    It will then remain to determine this cohomology via the
    second cochain complex in~(\ref{eqn:cacb}) above.

Fix $\g\in\Gamma$. 
    We will show that $\phi_{\g} d^* = d^* \phi_{\g}$ by
    checking that these two maps agree when evaluated
    on elements of the form $f\ot g$ where $f\in C^i(A,R)$
and $g\in C^j(B,R)$. 
Considering the $\N$-grading, 
the only elements of the resolution $P\ott Q$ 
    on which these maps
$(\phi_{\g} d^*) (f\ot g)$ and $(d^* \phi_{\g} )(f\ot g)$
    may be nonzero are sums of elements of the forms
    $(1\ot 1) \ot (a_1\ot\cdots \ot a_{i+1} ) \ot (b_1\ot \cdots \ot b_{j})$
    and  $(1\ot 1) \ot (a_1\ot\cdots \ot a_{i} ) \ot (b_1\ot \cdots
    \ot b_{j+1})$.
    (In~\cite{BW}, these two types of elements are
    labeled I and II.)
    We will check the results of evaluating on the second
    type of input; this is where the automorphism $\hg$ arises.
    The calculation for the first type of input is similar,
    but the automorphism does not arise since after applying $\phi_{\g}$,
    nonzero output occurs only where
    the differential is applied
    to the tensor factors in $A$ only.
Abusing notation, we will freely go back and forth
    between the expressions in~(\ref{eqn:pqab}) for elements of $P\ott Q$,
    and even further, between the expressions in~(\ref{eqn:homab})
    for cochains.
    In the calculation below, we will not record some of the output
    of the differential when it is clearly in
    the kernel of the next map to be applied, for example
$f(a\ot\cdots \ot a_{\ell} a_{\ell +1}\ot\cdots\ot a_i) = 0$
since $f$ has homological degree $i$ and this input has homological degree $i-1$. 
    Note that we may similarly assume without loss of generality that
    $a_1\ot\cdots\ot a_i \in (A^{\ot_R i})^{\g}$ since
    the $\Gamma$-degree of $f$ is $\g$,
 and $f$ takes values in $R$ which has degree~$1_{\Gamma}$. 
This is used in the first
    part of the calculation below, where
we also use the equality $\widehat{\g} (b_1) = t( a_1\ot\cdots\ot a_i , b_1) b_1$
(since $a_1\ot\cdots\ot a_i \in (A^{\ot_R i})^{\gamma}$): 
    \[
    \begin{aligned}
     & (d^*\phi_{\g}) (f\ot g)( (1\ot 1)\ot (a_1\ot\cdots\ot a_i)\ot
      b_1\ot\cdots\ot b_{j+1}))\\ & = 
      (\phi_{\g} (f\ot g)) (d((1\ot 1) \ot  (a_1\ot\cdots\ot a_i)\ot
      (b_1\ot\cdots\ot b_{j+1}))) \\
      & =  (\phi_{\g}(f\ot g)) \big( (-1)^i (1\ot \hg ( b_1))\ot
      (a_1\ot\cdots\ot a_i)
      \ot (b_2\ot \cdots b_{j+1})\\
      & \quad + (-1)^{i+j+1} \varepsilon(b_{j+1}) (1\ot 1)
         \ot (a_1\ot \cdots \ot a_i) \ot (b_1\ot\cdots\ot b_j)\\
      & \quad + (-1)^i \sum_{k=1}^j (-1)^k (1\ot 1)\ot (a_1\ot\cdots \ot a_i)
      \ot (b_1\ot\cdots\ot b_kb_{k+1}\ot\cdots\ot b_{j+1})\big) \\
      & = (-1)^{ij}(-1)^i \big(\varepsilon_B(\hg(b_1))
      f(a_1\ot\cdots\ot a_i)g(b_2\ot\cdots
      \ot b_{j+1}) \\
      & \quad + (-1)^{(i+1)(j+1)} \varepsilon(b_{j+1}) f(a_1\ot\cdots\ot a_i)
       g(b_1\ot\cdots \ot b_j) \\
      & \quad   + (-1)^{i(j+1)} \sum_{k=1}^j (-1)^k  
           f(a_1\ot\cdots\ot a_i) g(b_1\ot\cdots\ot
          b_kb_{k+1}\ot \cdots \ot b_{j+1} )\big) , 
      \end{aligned}
    \]
while taking into account the coefficients $R_{\widehat{\g}}$ in the expression
in~(\ref{eqn:cacb}) above, we have
      \[
      \begin{aligned}
        & (\phi_{\g} d^*) (f\ot g) ( (1\ot 1)\ot (a_1\ot\cdots\ot a_i)\ot
        b_1\ot\cdots\ot b_{j+1}))\\
        & = \phi_{\g}(d^*(f) \ot g +  f\ot d^*(g))
            ( (1\ot 1)\ot (a_1\ot\cdots\ot a_i)\ot
        b_1\ot\cdots\ot b_{j+1}))\\
        & = 
      (-1)^{i(j+1)} f(1\ot a_1\ot\cdots\ot a_{i})
        gd(1\ot b_1\ot\cdots\ot b_{j+1})\\
        & = (-1)^{i(j+1)}  f(a_1\ot\cdots\ot a_i)
        \varepsilon_B(\hg(b_1)) g(b_2\ot\cdots\ot b_{j+1}) \\
        & \quad + (-1)^{i(j+1)}(-1)^{j+1} \varepsilon(b_{j+1})
           f(a_1\ot\cdots\ot a_i)g(b_1\ot\cdots\ot b_j)\\
        & \quad + (-1)^{i(j+1)}  f( a_1\ot\cdots\ot a_i)
        \sum_{k=1}^{j} (-1)^k g(b_1\ot\cdots\ot b_kb_{k+1}
        \ot \cdots\ot b_{j+1}) , 
            \end{aligned}
      \]
      which agrees with the previous calculation.
We have proven the isomorphism~(\ref{eqn:cacb}). 

      Finally, we find the cohomology of the
      second cochain complex in~(\ref{eqn:cacb}).
      For each $\g\in\Gamma$,
      apply the K\"unneth Theorem~\cite[Theorem 3.6.3]{Weibel}  
to the tensor product of complexes
      $C(A,R)^{\g}$ and $C(B,R_{\hg})$.
      The hypotheses of the K\"unneth Theorem are satisfied since
      $R$ is hereditary and each
      $C(A,R)^{\g}$ is projective as an $R$-module. 
    Further, by hypothesis, each $\coh^i(A,R)$ is $R$-flat, 
and so the Tor term in the K\"unneth sequence 
(as in Remark~\ref{rk:Tor} below)
    is always 0.
    Therefore, as a graded $R$-module,  
the cohomology of the cochain complex~(\ref{eqn:cacb})
    is isomorphic to the tensor product of the cohomologies of its
    tensor factors.
    This is precisely the theorem statement.
\end{proof}

\begin{remark}\label{rk:Tor}
The hypothesis in
Theorem~\ref{thm:main} that
each $\coh^i(A,R)$ is flat over $R$ may be replaced by
  the hypothesis that $\coh^j(B,R_{\hg})$ is flat over $R$ for each $\g$, $j$.
If neither is flat over $R$, then
the isomorphism in the statement can be replaced by the following
short exact sequence of $R$-modules from the K\"unneth Theorem:
\[
\begin{aligned}
   0\rightarrow \bigoplus_{i+j=m, \gamma\in\Gamma}
  (\coh^i(A,R)^{\gamma}\ot_R \coh^j(B,R_{\hg})) \rightarrow
 &  \coh^m(\ttpt,R) \rightarrow \\
 & \bigoplus_{i+j=m-1, \gamma\in \Gamma}
   \Tor^R_1(\coh^i(A,R)^{\gamma}, \coh^j(B,R_{\hg})) \rightarrow 0 .
\end{aligned}
\]
\end{remark}

The following corollary should be compared with~\cite[Theorem 3.7]{BO}, 
which is for graded modules. 
Under our additional hypothesis in the corollary, that is that $R$ is also graded
as a $B$-module, there are $B$-module isomorphisms $R_{\hg}\cong R$
for all $\gamma\in\Gamma$ since $B^{\lambda}$ is in $\Ker(\varepsilon_B)$
for all $\lambda\neq 1_{\Lambda}$ and 
the automorphism $\hg$ acts as the identity on 
$B^{1_{\Lambda}}$.
In this case, the proof techniques can be pushed yet further to yield 
the multiplicative result Theorem~\ref{thm:cup} below, 
just as in~\cite[Theorem 3.7]{BO}. 

\begin{cor}\label{cor:main}
Let $R$ be a hereditary commutative ring that is both a
$\Gamma$-graded $A$-module and a $\Lambda$-graded $B$ module,
of degrees $1_{\Gamma}$ and $1_{\Lambda}$, respectively.
Assume that for each $i\geq 0$, $\coh^i(A,R)$ is flat over $R$.
Then for all $m\geq 0$, there is an isomorphism of $R$-modules:
\[
   \coh^m(\ttpt, R)\cong   
   \bigoplus_{i+j = m}
    (\coh^i(A,R) \ot_R
    \coh^j(B,R)).
\]                  
\end{cor}

\begin{proof}
Since $R$ is a $\Lambda$-graded $B$-module, there is an induced
$\Lambda$-grading on cohomology via the bar resolution.
The statement here is just the corresponding version
of the statement in Theorem~\ref{thm:main}, using the
isomorphisms $R_{\hg}\cong R$ described above for all $\gamma\in\Gamma$
in this case that $R$ is also $\Lambda$-graded. 
\end{proof}

\section{Cup products in the graded module setting}\label{sec:cup}

In this section, we describe the multiplicative structure of $\coh(A\ott B, R)$
in the case that $R$ is graded both as an $A$-module and as a $B$-module
as in Corollary~\ref{cor:main}. 
In the next section, we handle a class of examples where this does not hold,
namely the Taft orders.

The proof of~\cite[Theorem 3.7]{BO} adapts with no difficulties to this setting
to establish Theorem~\ref{thm:cup} below.
The proof should also be compared with~\cite[Theorem 4.7]{BW} for 
Hochschild cohomology.  
We give a brief summary for completeness, 
expressing the cup product 
at the cochain level via the isomorphism~(\ref{eqn:cacb}) in the proof of 
Theorem~\ref{thm:main}.
Let $f\in C^i(A,R)^{\g}$, $f'\in C^{i'}(A,R)^{\g '}$,
$g\in C^j(B, R)^{\lambda}$, $g'\in C^{j'} (B, R)^{\lambda'}$.
In particular, the $\N\times\Gamma$-degree of $f$ is $(i,\gamma)$,
that of $g$ is $(j,\lambda)$,
and similarly for $f'$, $g'$. 
Set
\begin{equation}\label{eqn:cup}
  \begin{aligned}
&  (f\ot g)\smile (f'\ot g') (a_1\ot\cdots \ot a_{i+i'}\ot b_1\ot
  \cdots \ot b_{j+j'}) \\ & =
  (-1)^{ii' +ij' + jj'}
     t^{-1}( a_{i+1}\ot \cdots \ot a_{i+i'},  b_1\ot\cdots\ot b_j) \cdot \\
 & \hspace{1cm}f(a_1\ot \cdots \ot a_i)
      g(b_1\ot\cdots\ot b_j) 
      f'(a_{i+1}\ot\cdots\ot a_{i+i'})
      g'(b_{j+1}\ot\cdots\ot b_{j+j'}) ,
  \end{aligned}
\end{equation}
for homogeneous elements $a_1,\ldots, a_{i+i'}$ in $A$
and $b_1,\ldots, b_{j+j'}$ in $B$. 

For the following statement, the bicharacter twisting map $\widetilde{t}$
is defined by
\[
    \widetilde{t}((i,\gamma), (j, \lambda)) =
    (-1)^{ij} t^{-1}(\gamma,\lambda) 
\]
for all $i,j\in \N$, $\gamma\in\Gamma$, $\lambda\in\Lambda$. 

\begin{theorem}\label{thm:cup}
Let $R$ be a hereditary commutative ring
that is both a $\Gamma$-graded $A$-module and a $\Lambda$-graded $B$-module,
of degrees~$1_{\Gamma}$ and $1_{\Lambda}$, respectively.
Assume that for each $i\geq 0$, 
  $\coh^i(A,R)$ is flat over $R$.
There is an isomorphism of $\N$-graded rings  
  \[
  \coh (A\ott B , R)\cong 
    \coh (A, R) \ot^{\widetilde{t}}_R \coh (B, R) .  
   \]
\end{theorem}

\begin{proof}
In the notation above, set $F = f\ot g$ and $F'=f'\ot g'$. 
We will show that the cup product $F\smile F'$
at the chain level as a function on the twisted tensor product
$P\ott Q$ of two bar resolutions agrees with~(\ref{eqn:cup}).
Let $\Delta: P\ot^t\! Q\rightarrow  (P\ott Q)\ot_R (P\ot^t\! Q)$
be the composition
\[
   P\ot^t\! Q \stackrel{\Delta_A\ot \Delta_B}
     {\relbar\joinrel\relbar\joinrel\relbar\joinrel\relbar\joinrel\longrightarrow} 
   P\ot_R P \ot_R Q\ot_R Q \stackrel{1\ot\tilde{\tau}\ot 1}
  {\relbar\joinrel\relbar\joinrel\relbar\joinrel\longrightarrow}
   P\ot_R Q \ot_R P\ot_R Q
\]
where $\Delta_A$, $\Delta_B$ are standard diagonal maps on
bar resolutions, and $\tilde{\tau}: P\ot_R Q\rightarrow Q\ot_R P$ is the
standard twisting map on these resolutions given by iterating $\tau^{-1}$.
Then by construction, $\Delta$ is a chain map lifting the canonical isomorphism 
$R\stackrel{\sim}{\longrightarrow} R\ot_R R$.
We set $F\smile F' =  (F\ot_R F')\circ \Delta$.
The claimed formula follows.
\end{proof}

\begin{remark}
A similar formula to~(\ref{eqn:cup}) above applies in case the bar resolutions
on $A$ and $B$ are replaced by arbitrary projective resolutions $P$ and $Q$,
where $P$ is compatible with the twisting, given diagonal maps $\Delta_P$
and $\Delta_Q$. 
\end{remark}

\begin{cor}\label{cor:fg}
Let $R$ be a hereditary commutative ring
that is both a $\Gamma$-graded $A$-module and a $\Lambda$-graded $B$-module,
of degrees~$1_{\Gamma}$ and $1_{\Lambda}$, respectively.
Assume that for each $i\geq 0$,
  $\coh^i(A,R)$ is flat over $R$. 
If $\coh(A,R)$ and $\coh(B,R)$ are both finitely generated as $R$-orders,
then $\coh(\ttpt, R)$ is a finitely generated $R$-order.
\end{cor}

\begin{proof}
The result follows from Theorem~\ref{thm:cup} by
taking the generating set to consist of tensor products of elements
of generating sets of $\coh(A,R)$ and of $\coh(B,R)$
with the respective identity elements. 
\end{proof}

\begin{example}
(Some quantum complete intersections over algebraic integers.)
Let $\omega$ be an $m$th root of unity in $\C$ and $R = \Z[\omega]$.
Let $S$ be the augmented $R$-order given by generators $x_1,\ldots,x_n$, 
relations
\[
    x_ix_j = \omega_{ij} x_jx_i \ \ \ \mbox{ and } \ \ \
    x_i^{N_i} = 0
\]
for all $i < j$ where $N_i\geq 2$ and each $\omega_{ij}$ ($i < j$) is a power of $\omega$, 
and augmentation $\varepsilon(x_i) = 0$ for all $i$.
Note that $S$ is isomorphic to an iterated bicharacter twisted tensor
product of its suborders $R[x_i]/(x_i^{N_i})$, each
graded by $\Gamma = \Z/m\Z$ with $x_i$ in degree 1.
The values of the bicharacter $t$ on the grading group $\Gamma$ 
are given by powers of $\omega$. 
Further, $R$ is a graded $S$-module
concentrated in degree 0 as an element of the grading group $\Z/m\Z$. 
The standard minimal projective $R[x_i]/(x_i^{N_i})$-module resolution of $R$ 
(see~(\ref{eqn:P-std}) below with $n=N_i$, $x = x_i$) 
gives $\coh^{\ell}(R[x_i]/(x_i^{N_i}),R)\cong R$ for all $\ell \geq 0$.
By Theorem~\ref{thm:cup}, $\coh(S,R)$ is isomorphic to a twisted
tensor product of all $\coh(R[x_i]/(x_i^{N_i}), R)$, a straightforward
generalization of~\cite[Theorem 4.1]{MPSW} and of a special case
of~\cite[Theorem 5.3]{BO}:
When $N_i >2$ for all $i$,
\[
  \coh (S,R) \cong R\langle\xi_i, \eta_i \rangle /
   (\xi_i\xi_j - \omega_{ji}\xi_j\xi_i , \ 
    \eta_i\xi_j - \omega_{ji}^{N_j} \xi_j\eta_i , \
   \eta_i \eta_j + \omega_{ji}\eta_j\eta_i )
\]
where $\deg \xi_i = 2$ and $\deg \eta_i = 1$.
In particular, $\eta_i^2 = 0$ for each $i$.
For any $i$ for which $N_i = 2$, the cocycle $\eta_i^2$ is a nonzero scalar multiple of
$\xi_i$ and so $\xi_i$ may be removed from the generating set with no loss. 
\end{example}

In contrast with the above quantum complete intersection example,
in other settings such as the Taft orders in
the next section,
cohomology can have torsion elements and its ring structure
need not be a twisted tensor product of cohomology of constituent rings.


\section{Taft orders over cyclotomic integers}\label{sec:Taft}

In this section we handle Taft orders,
to which the results of Section~\ref{sec:cup} do not apply.
Instead, we work directly with cochains on twisted tensor product resolutions, 
illustrating how the twisted tensor product methods
of Sections~\ref{sec:ttp} and~\ref{sec:additive}
may be pushed further to obtain the multiplicative structure on cohomology
for additional classes of examples.

Let $n$ be an integer ($n\geq 2$), 
let $\om$ be a primitive $n$th root of unity in $\C$,
and $R = \Z[\om]$.
Then $R$ is a Dedekind domain, and in particular is hereditary.
Let $T_n$ be the augmented $R$-order given by generators $x$ and $g$,
relations 
\[
   x^n=0, \ \ g^n=1, \ \ gx=\om xg ,
\]
and augmentation given by $\varepsilon(x)=0, \ \varepsilon(g) =1$.
As a free $R$-module, $T_n$ has basis
$\{x^ig^j\mid 0 \leq i,j\leq n-1\}$. 
As a Hopf order, the coproduct $\Delta$ on $T_n$ is determined by
\[
   \Delta(g) = g\ot g , \ \ \ \Delta(x) = x\ot 1 + g\ot x . 
\]
We will find the cohomology
$
   \coh(T_n , R) = \Ext^*_{T_n}(R,R) .
$

Note that $T_n$ is isomorphic to a twisted tensor product
$\ttp$ of its two suborders 
\[
   A = R[x]/(x^n) \ \ \ \mbox{ and } \ \ \  
B = R G,
\]
where $G = \langle g\rangle\cong \Z/n\Z$: 
The twisting map $\tau: B\ot_R A \rightarrow A\ot_R B$ is given by 
\[
  \tau(g\ot x) = \om x\ot g , 
\]
extended to $B\ot_R A$ by requiring compatibility~(\ref{eqn:twisting-mult}) 
with multiplication on $A$ and $B$. 
This can be viewed as bicharacter twisting: 
Let $\Gamma = \Lambda = G$,
which grades $B = RG$ canonically and grades $A$ by 
giving $x$ the degree $g$.
The bicharacter $t$ is defined by
  $
  t(g^i, g^j ) = \omega^{ij} .
  $
The $T_n$-module $R$, restricted to $A$, is
compatible with $\tau$ via the canonical isomorphism
$\tau_{B,R} : B\ot_R R \stackrel{\sim}{\longrightarrow} R\ot_R B$.

\begin{theorem}\label{thm:Taft}
  The cohomology of the Taft $R$-order $T_n$ 
in degree $m$ with coefficients in $R =\Z[\omega]$ is
  \[
  \coh^m(T_n,R) \cong \left\{
  \begin{array}{ll}
    R\oplus  \left( R/(\omega -1) \right)^{\oplus \frac{m}{2}}
    \oplus \left( R/(n) \right)^{\oplus \frac{m}{2}} ,
    & \mbox{ if $m$ is even}\\
    0 , & \mbox{ if $m$ is odd} .
  \end{array}\right.
  \]
\end{theorem}

\begin{proof}
This is a direct application of Theorem~\ref{thm:main}.
We will use the following resolutions to find the needed cohomology
$\coh^i(A,R)^{\gamma}$ and $\coh^j(B, R_{\widehat{\gamma}})$.

Let $P$ be the projective resolution of $R$ as an $A$-module given by:
\begin{equation}\label{eqn:P-std}
\xymatrix{
P: \quad \quad \cdots \ar[r] &
   A\ar[r]^{x^{n-1}\cdot} &
   A\ar[r]^{x\cdot} &
   A\ar[r]^{x^{n-1}\cdot} &
   A\ar[r]^{x\cdot} &
    A \ar[r]^{\varepsilon} & R \ar[r] & 0 .
}
\end{equation}
Since $\varepsilon(x)=0$, the Hom complex
$\Hom_A (P, R)$ is, in each degree $i$, $\Hom_A(A,R)\cong R$,
with differentials all 0.
So $\coh^i(A,R) \cong R$ for each $i$.
In particular, $\coh ^i(A,R)$ is flat over $R$. 
Note that $P$ is compatible with $\tau$ via the maps
$\tau_{B,P_i}: B\ot_R P_i \rightarrow P_i\ot_R B$ given by
\[
   \tau_{B, P_i}( g^{\ell}\ot  x^j) =  \left\{
    \begin{array}{rl}
     \omega^{j\ell}x^j\ot g^{\ell} , & \mbox{ if $i$ is even}\\
    \omega^{(j+1)\ell} x^j \ot g^{\ell} , & \mbox{ if $i$ is odd} .
   \end{array}\right.
\]

Let $Q$ be the projective resolution of $R$ as a $B$-module given by:
\[
\xymatrix{
 Q: \quad \quad \cdots \ar[r]
    & B\ar[r]^{v\cdot} & 
    B\ar[r]^{u\cdot} &
    B\ar[r]^{v\cdot} &
    B\ar[r]^{u\cdot} & 
     B\ar[r]^{\varepsilon} & R \ar[r] & 0   , 
}
\]
where $u=g-1$ and $v=1+g+g^2+\cdots + g^{n-1}$. 
The $T_n$-module structure on each $P_i\ot_R Q_j$
arising from the twisting maps $\tau_{B, P_i}$ defined above is as follows. 
The action of $A$ is on the leftmost tensor factor $P_i$ only. 
The action of $B$ involves the twisting maps $\tau_{B, P_i}$:
For all $h\in G$ and $0\leq j\leq n-1$,
\begin{equation}\label{eqn:action}
  g^{\ell}\cdot (x^j\ot h) =
     \left\{\begin{array}{rl}
     \omega^{j\ell} x^j\ot g^{\ell}h, & \mbox{ if $i$ is even}\\
    \om^{(j+1)\ell} x^j\ot g^{\ell}h , & \mbox{ if $i$ is odd} .
   \end{array}\right.
\end{equation}
In this way, $P\ot_R Q$ is a complex of $T_n$-modules.

Finally, we show that each $P_i\ot_R Q_j$ is 
a projective $T_n$-module (in fact, free of rank one):
Fix $i,j$. Note first that 
$1_A\ot 1_B$ generates $P_i\ot_R Q_j$ as a $T_n$-module. 
Define a $T_n$-module homomorphism $f$ from
$T_n$ to $P_i\ot_R Q_j$ by sending $1_{T_n}$ to $1_A\ot 1_B$;
since $T_n$ is free as a $T_n$-module, 
this extends to a $T_n$-module homomorphism as desired.
Considering the above actions, $f$ is invertible, and so $P_i\ot_R Q_j$
is a free $\ttpt$-module of rank one. 
By Theorem~\ref{thm:ttpr}, the total complex of $P\ot_R Q$ is a projective resolution
of $R$ as a $T_n$-module, denoted $P\ott Q$. 

Next note that the hypotheses of Theorem~\ref{thm:main} hold.
  In particular, $R$ is a graded $A$-module, and 
as explained above, for each $i$, the cohomology
  $\coh ^i (A,R)$ is flat over $R$.
  As a $B$-module,  $R$ is {\em not} graded, since for example
  $g$ acts as the identity map on $R$. 
  It can be checked, using the above resolutions, that
  \begin{eqnarray*}
 \coh^i(A,R)^{g^s} & \cong & \left\{\begin{array}{rl}
    R , & \mbox{ if $i$ is even, } s = 0\\
    R , & \mbox{ if $i$ is odd, } s= -1 \\
    0 , & \mbox{ otherwise}
    \end{array}\right. \\
 \coh^j(B,R_{\widehat{1}}) & \cong & \left\{\begin{array}{rl}
    R , & \mbox{ if } j = 0\\
    R/(n) , & \mbox{ if $j$ is even, } j>0 \\
    0 , & \mbox{ if $j$ is odd}
    \end{array}\right. \\
 \coh^j(B, R_{\widehat{g^{-1}}}) & \cong & \left\{\begin{array}{rl}
    0, & \mbox{ if $j$ is even}\\
    R/(\omega -1) , & \mbox{ if $j$ is odd} 
    \end{array}\right.
    \end{eqnarray*}
Applying Theorem~\ref{thm:main} now, we see from
the above isomorphisms that 
    the cohomology is concentrated in even degrees.
    For $m$ even, 
    the identity component of $\coh^m(A\ott B, R)$
    is the sum over all even $i$ and $j$ for which $i+j =m$ of 
    \[
    \coh^i(A,R)^1 \ot_R \coh^j(B,R_{\widehat{1}}) ,
    \]
    which is $R\oplus \frac{R}{(n)}^{\oplus \frac{m}{2}}$.
The $g^{-1}$ component
    is the sum over all odd $i$ and $j$ for which $i+j=m$ of
    \[
    \coh^i(A,R)^{g^{-1}}\otimes_R \coh^j(B, R_{\widehat{g^{-1}}}) ,
    \]
    which is $\frac{R}{(\omega - 1)} ^{\oplus \frac{m}{2}}$.
The $g^s$ components for $s\neq 0, \ -1$ are all zero.
Thus $\coh^m (\ttpt , R)$ is as claimed.
\end{proof}

Next we wish to describe cup products. 
Theorem~\ref{thm:cup} does not apply since $R$ is not graded as a $B$-module.
We will instead 
give an explicit description of $\coh(T_n, R)$ as a ring
by looking more closely at the projective resolution $P\ott Q$
from the proof of Theorem~\ref{thm:Taft}.

\begin{theorem}\label{thm:Taft-cup}
For all $n>2$, the cohomology of the Taft $R$-order $T_n$
with coefficients in $R = \Z[\omega]$ is 
\[
      \coh(T_n,R)\cong R [ \alpha,\ \beta,\ \gamma ] /
     \big( \gamma^2 , \ (\omega -1)\gamma , \ n\beta \big) ,
\]
and for $n=2$ it is 
\[
  \coh(T_2,R)\cong R [ \alpha, \ \beta, \ \gamma ] /
     \big( \gamma^2 - \alpha\beta, \ 2\gamma, \ 2\beta \big) ,
\] 
where $\alpha, \ \beta, \ \gamma$ are generators in degree~2.
\end{theorem} 

\begin{proof}
We will first derive an explicit formula for the differential on
$P\ott Q$. 
For notational clarity, let $\phi_{ij} = 1_A\ot 1_B$ in $P_i\ott Q_j$.
Then we may write
\[
  P_i\ot_R Q_j = T_n \phi_{ij} ,
\]
the free $T_n$-module of rank one with free basis element $\phi_{ij}$. 
Similarly, write
$P_i=A\phi_i$ and $Q_j=B\psi_j$ where
$\phi_i$ and $\psi_j$ denote free basis elements, 
and identify $\phi_{ij}$ with $\phi_i \ot \psi_j$. 
To obtain a formula for the differential,
we use the action~(\ref{eqn:action})
to write all expressions in terms of the free basis elements $\phi_{ij}$:
\begin{eqnarray*}
d(\phi_{ij})&=& d(\phi_{i}\ot\psi_j) \\
&=& d(\phi_{i})\ot\psi_j + (-1)^{i} \phi_{i}\ot d(\psi_j)\\
  &=&    \left\{  \begin{array}{ll} 
     x\phi_{i-1}\ot \psi_j - \phi_i\ot (g-1)\psi_{j-1} , &
       \mbox{ if $i$ odd, $j$ odd }\\
     x\phi_{i-1}\ot\psi_j - \phi_i\ot (1+g+g^2+\cdots g^{n-1})\psi_{j-1} ,
   & \mbox{ if $i$ odd, $j$ even } \\
  x^{n-1}\phi_{i-1}\ot \psi_j + \phi_i\ot (g-1)\psi_{j-1} ,
   & \mbox{ if $i$ even, $j$ odd } \\
   x^{n-1}\phi_{i-1}\ot\psi_j + \phi_i\ot (1+g+g^2+\cdots +g^{n-1})
   \psi_{j-1} , & \mbox{ if $i$ even, $j$ even } 
  \end{array}\right.  \\
  & & \\ & = &
   \left\{\begin{array}{ll}
   x\phi_{i-1,j} - (\om^{-1}g -1) \phi_{i,j-1}, &
   \mbox{ if $i$ odd, $j$ odd } \\
  x\phi_{i-1,j} - (1+\om^{-1}g + \om^{-2}g^2 + \cdots + \om^{1-n}g^{n-1})
   \phi_{i,j-1}, &\mbox{ if $i$ odd, $j$ even }\\
  x^{n-1}\phi_{i-1,j} + (g-1)\phi_{i,j-1} , 
   & \mbox{ if $i$ even, $j$ odd }\\
  x^{n-1}\phi_{i-1,j} + (1+g+g^2+\cdots + g^{n-1})\phi_{i,j-1},
  &\mbox{ if $i$ even, $j$ even }.
  \end{array}\right.
\end{eqnarray*}

\quad

For use in calculations below, note that since $\om$ is a primitive $n$th 
root of unity, by factoring the polynomial 
$1+z+z^2 + \cdots + z^{n-1}$ and substituting $z=g$, $z=\om^{-1} g$, we obtain
\begin{eqnarray*}
   1 + g + g^2 + \cdots + g^{n-1} &=& (1-\om g) (1 - \om^2 g) \cdots (1 - \om^{n-1}g), \\
   1 + \om^{-1}g + \om^{-2}g^2 + \cdots + \om^{1-n}g^{n-1} &=& (1-g)(1-\om g)\cdots
     (1 - \om^{n-2}g) .
\end{eqnarray*}

Apply $\Hom_{T_n}( - , R)$ to each $P_i\ot_R Q_j$
and note that for each $i,j$, 
\[
  \Hom_{T_n}(P_i\ot_R Q_j , R)\cong
  \Hom_{T_n}(T_n \phi_{ij}, R)\cong \Hom_R (R\phi_{ij},R) ,
  \]
that is the free $R$-module of rank one dual to $R\phi_{ij}$. 
Denote by $\phi_{ij}^*$ (in $\Hom_R (R\phi_{ij},R)$) 
the dual basis element to $\phi_{ij}$.

Using the above formula,
the differential on the total complex of
$\Hom_{T_n}(P_i\ot_R Q_j , R)$ is given by:
\[
  d^*(\phi^*_{ij}) = \left\{\begin{array}{ll}
    (1-\om^{-1})\phi^*_{i,j+1} , & \mbox{ if $i$ odd, $j$ even}\\
   n \phi_{i,j+1}^*, & \mbox{ if $i$ even, $j$ odd}\\
   0 , & \mbox{ if $i+j$  even}
  \end{array}\right.
\]

\quad

For clarity, we verify that the cohomology $\coh(T_n,R)$ 
obtained explicitly from this complex agrees with
Theorem~\ref{thm:Taft}:
Since $\omega (1-\omega^{-1})= \omega -1$,
the elements $1-\omega^{-1}$ and $\omega -1$
generate the same ideal of $R$. 
In case $n=2$,
the ideals $(2)$ and $(\omega -1)$ are equal in $R = \Z$,
so the statement simplifies in that case.
More generally, there are other relations between the ideals
$(n)$ and $(\omega -1)$ in $R$, for example if $n=3$ then $((\omega-1)^2)=(3)$. 
If $m$ is odd, the cohomology in degree $m$ is 0:
it is clear
 from the above formula for the differential $d^*$ on the cochain complex
  that the kernel is 0 in odd degrees.
  Now suppose $m$ is even. Note that $d^*(\phi^*_{ij}) = 0$ whenever $i+j=m$
  by the above formula for $d^*$.
  Taking the quotient by the image of the differential from degree $m-1$,
  the formula for $d^*$ above yields
  \begin{equation}\label{eqn:coh-revisited}
  R\phi^*_{m0}  \oplus
  \left(\bigoplus_{i+j=m, \ i \mbox{ odd}} \frac{R}{(\omega -1)}\phi^*_{ij} \right)
  \oplus \left(\bigoplus_{i+j=m, \ i \mbox{ even}, \ j>0}
          \frac{R}{(n)}\phi^*_{ij}\right) ,
  \end{equation}
          a direct sum of one copy of $R$ with $\frac{m}{2}$ copies
          of $R/(\omega-1)$ and $\frac{m}{2}$ copies of $R/(n)$,
in accordance with Theorem~\ref{thm:Taft}.

Next we consider the graded commutative ring structure
on $\coh(T_n,R)$.
We will find it by composing chain maps corresponding to
cocycles in degree 2;
our computations will verify that indeed all
generators are in degree 2.
We will need some more notation:
Let $C_{\bu}$ be the total complex $P\ott Q$, so that
\[
   C_m = \bigoplus_{i+j=m} T_n \phi_{ij}.
\]
By the formula for differentials above, we see that
$\phi^*_{20}, \phi^*_{11}, \phi^*_{02}$ are cocycles in degree~2.
Let $f_{20}, f_{11}, f_{02}$ denote choices of corresponding
chain maps from $C$ to $C$ of degree $-2$.
For example, $f_{20}$ fits into a commutative diagram: 
\[
\xymatrix{
  \cdots \ar[r] & C_4\ar[r]^{d}\ar[d]_{f_{20}} & C_3\ar[r]^{d}\ar[d]_{f_{20}}
    & C_2\ar[r]^{d}\ar[d]_{f_{20}} \ar[dr]^{\hspace{-.1cm}\phi_{20}^*} 
    & C_1\ar[r]^{d} 
    & C_0\ar[r]^{\varepsilon} &   R \ar[r] & 0 \\
  \cdots \ar[r]^{d} & C_2\ar[r]^{d} & C_1\ar[r]^{d} & C_0 \ar[r]^{\varepsilon} 
  & R\ar[r] & 0 & &   
}
\]
Calculations show that we may choose
\[
   f_{20}(\phi_{ij}) = \phi_{i-2 , j} \ \ \mbox{ for all }i,j
\] 
(where a negative index indicates the value is zero). 
Similarly, we may take   
\[ 
  f_{02}(\phi_{i,j}) = \phi_{i,j-2} \ \ \mbox{ for all }i,j.
\]
To find a suitable chain map $f_{11}$,
we first define the following elements $a,b$ of $T_n$:
\begin{eqnarray*}
   a & := & (1-\om g) (1-\om ^2 g)\cdots (1-\om^{n-2} g), \\
   b & := & (1-\om ^2 g)(1-\om^3 g)\cdots (1-\om^{n-1}g) x^{n-2},
\end{eqnarray*} 
where if $n=2$, we interpret this to mean $a:=1$, $b:=1$. 
Calculations show that we may take  
\[
   f_{11}(\phi_{ij}) = \left\{\begin{array}{ll}
   b\phi_{i-1, j-1}, & \mbox{ if $i$ even, $j$ even}\\
   x^{n-2}\phi_{i-1,j-1}, & \mbox{ if $i$ even, $j$ odd}\\
  a \phi_{i-1,j-1},& \mbox{ if $i$ odd, $j$ even}\\
   \phi_{i-1,j-1},& \mbox{ if $i$ odd, $j$ odd} .
  \end{array} \right.
\]
Note that when $n=2$, this simplifies to $f_{11}(\phi_{ij}) = \phi_{i-1,j-1}$.
By composing these maps at the chain level, 
we find that 
\[
  \phi^*_{ij} \phi^*_{02} = \phi^*_{i, j+2}
  \ \ \mbox{ and } \ \ 
  \phi^*_{ij} \phi^*_{20} = \phi^*_{i+2, j}.
\]
If $n>2$, since $ax^{n-2}$ and $b$ are both in the kernel of $\varepsilon$, 
we have $(\phi^*_{11})^2 =0$,
while if $n=2$, we have $(\phi_{11}^*)^2 = \phi^*_{20}\phi^*_{02}$.
Now denote by $\alpha, \beta, \gamma$ the degree 2 elements of cohomology  
corresponding to $\phi_{20}^*$,
$\phi_{02}^*$, $\phi_{11}^*$ at the chain level, respectively. 
The above relations among $\alpha, \beta, \gamma$
together with the description~(\ref{eqn:coh-revisited}) of the cohomology
in each even degree~$m$, including the torsion details, now complete the proof.
\end{proof}

\end{document}